\documentclass[12pt]{article}
\usepackage{graphicx}
\usepackage{amsmath}
\usepackage{amssymb}
\usepackage{vmargin}
\usepackage{amsthm}

\usepackage[mathscr]{eucal}
\usepackage{amsbsy}
\usepackage{amsmath}
\usepackage{graphics}
\usepackage{epsfig}

\usepackage{latexsym}

\newtheorem{Theorem}{Theorem}

 \newtheorem{Corollary}[Theorem]{Corollary}
\newtheorem{Lemma}[Theorem]{Lemma}
\newtheorem{Proposition}[Theorem]{Proposition}
\newtheorem{Remark}[Theorem]{Remark}

\def\k{\varkappa}

\def\=A8{\"o}

\title{{\bf  Point Process Approach to the Winner Problem} }
\author{{Youri \sc Davydov${}^{1}$}\quad  and\quad { Vladimir \sc Rotar ${}^{2}$}}
\date
{\footnotesize  ${}^{1}$  Laboratoire Paul Painlev\'e \\
 Universit\'e de Lille, France,\\
 and Faculty of Mathematics and Computer Sciences\\
 St. Petersburg State University, Russia\\
 Email: youri.davydov@univ-lille.fr\\
\vspace{3pt}
${}^{2}$ \ Department of Mathematics \\
 University of California at San Diego, USA and \\
 the National University, USA \\
  Email: vrotar@math.ucsd.edu }

\begin{document}

\maketitle

\begin{quote}
{\bf Abstract.}
We consider a limit theorem for a triangular array of point processes generated by non-identically distributed random variables, and apply the result  for the analysis of the limiting behavior of the Argmaximum of independent random variables, as well as for some step processes.
\end{quote}


\setcounter{Theorem}{0}

\section{Introduction}\label{secintro}

\renewcommand{\theequation}{\thesection.\arabic{equation}}
\setcounter{equation}{0}

 In our previous paper \cite{d-r}, we considered the asymptotic behavior of the \textit{Argmaximum} of a large number of  independent random variables (r.v.'s). In \cite{d-r}, we  called this problem \textit{the winner problem}.

 If the r.v.'s are identically distributed,  the answer is obvious: from the very beginning, the distribution of the Argmaximum is uniform, while for non-identically r.v.'s, the problem turned out to be non-trivial.

 Clearly, for the result to be substantial, the tail distributions  should have similar in a sense character.  Since the distribution of an Argmaximum is invariant under a strictly increasing transformation of the r.v.'s, we may suppose from the very beginning that the mentioned tails satisfy the condition of regular variation. Since the last condition is strongly connected with the convergence of empiric point processes (see, for example Resnick \cite{r}), it could be  natural to use this connection for the analysis of the limiting behavior of the Argmaximum distribution.

 To this end, firstly we prove a theorem (Theorem \ref{th1}) on the convergence of point processes for a triangular array of non-identically distributed r.v.'s, This theorem generalizes the well known Resnick's theorem (see, Proposition 3.21 in \cite{r}).

 Next, we derive from Theorem \ref{th1} a number of corollaries on the behavior of the Maximum, Argmaximum, and step processes.
 Our previous result about
 argmax is also partially covered by this theorem.

\section{Point Processes}\label{secpointprocesses}

\renewcommand{\theequation}{\thesection.\arabic{equation}}
\setcounter{equation}{0}

In this section, we state a limit theorem for empiric point processes generated by a triangular array of independent r.v.'s.


Let   $\mathbf{X}=\{X_{n,j},\;j\leq n;\;n\in \mathbb{N}\}$  be a triangular array of independent non-negative r.v.'s, and let  $\mathcal{P}_{n,j}$ be the distribution of r.v.'s $X_{n,j}$. We assume that all $\mathcal{P}_{n,j}$  are non-atomic.

Denote by   $\mathbf{C}=\{c_{n,j},\;j\leq n;\;n\in \mathbb{N}\}$  a triangular array of positive constants for which the following set of conditions is true:

$\mathbf{H_1}:$  {\it As  $n\rightarrow \infty$
\begin{equation}\label{c}
d_n:= \sum_{j=1}^n c_{n,j}  \rightarrow \infty;\;\;\;\;\;\;
\frac{\max_{j\leq n}\{c_{n,j}\}}{d_n} \rightarrow 0,
\end{equation}
and measures on  $[0,1]$
$$
\mu_n:= \sum_{j=1}^n 	\frac{c_{n,j}}{d_n}\delta_{\{\frac{j}{n}\}}
$$
 weakly converge to a measure $\mu$:
\begin{equation}\label{h1}
	\mu_n{\Longrightarrow} \mu.
\end{equation}
}
(Here,  $\delta_{\{a\}}$ is a measure concentrated at a point $a$.)

\begin{Remark}
 In virtue of (\ref{c}),   measure $\mu$ is non-atomic. Therefore the
	 convergence (\ref{h1}) is equivalent to the following: for any $t\in [0,1]$, weakly,
	$$
\frac{1}{d_n} \sum_{j\leq nt} c_{n,j} \rightarrow	\mu([0,t]).
	$$
\end{Remark}	

\vspace{5pt}

Now on, we will use notation $\mathbb{R}_+ = (0,\infty]$ for a half-line with an added point $+\infty$ and provided with a metric $\rho(x,y)= |\frac{1}{x}-\frac{1}{y}|,\;x,y >0.$

In this case, a subset $A \subset R_+$ is compact if it is closed in the usual topology and separated from zero. The point $+\infty$ has been added for convenience in order that $(R_+, \rho)$ will become a complete metric space.

All measures we will consider will be finite on $[t,\infty]$ for any $t>0$ and will not have an atom at the point $+\infty$.

Denote by $\mathbb{C}_{K,t}^+$   the set of continuous in $\rho$ functions with a compact support in $[t, \infty)$. Set  $\mathbb{C}_{K}^+  =
\cup_{t>0}\mathbb{C}_{K,t}^+.$

Let us recall that measures $\tau_n$  on   $\mathbb{R}_+$     converge to a measure   $\tau $      in a     $vague$-topology if for any $f\in \mathbb{C}_{K}^+$
 \begin{equation}\label{vague1}
    \int f\,d\tau_n \rightarrow \int f\,d\tau.
 \end{equation}

In our case, in essence, this means that for any $t>0$, the restrictions of measures $\tau_n$ on [$t,\infty)$ converges in the usual topology to the corresponding  restriction of $\tau$. If measure $\tau$ is non-atomic, then (\ref{vague1}) is equivalent to
\begin{equation}\label{vague2}
\tau_n([t,\infty)\to\tau([t,\infty) \,\,\,\, \text{as}\,\,\,\, n\to\infty\,\,\, \text{for any}\,\,\,\, t>0.
\end{equation}
We suppose that collections   $\mathbf{X}$ and  $\mathbf{C}$  satisfy the following {\bf regular variation} condition:

$\mathbf{H_2}:$  {\it For a non-atomic  measure   $\gamma$  on $\mathbb{R}_+$, for any  $j$, as  $n\rightarrow \infty$
\begin{equation}\label{h2.1}
\frac{d_n}{c_{n,j}}\mathcal{P}_{n,j} \stackrel{vague}{\longrightarrow} \gamma,
\end{equation}
and this convergence is uniform in a sense that for any  $t>0$,  as $n\rightarrow \infty$,
\begin{equation}\label{h2.2}
\Delta_n(t):=  \sup_{f\in \mathbb{C}_{K,t}^+}\max_{j\leq n}
\left|\frac{d_n}{c_{n,j}}\int f\,d\mathcal{P}_{n,j} - \int f\,d\gamma
\right| \rightarrow 0.
\end{equation}
}
On  $[0,1]\times\mathbb{R}_+$, we define point processes

\begin{equation}
	\zeta_{n} = \sum_{j=1}^n \delta_{\{(\frac{j}{n}, \;X_{n,j})\}}.
\end{equation}

\begin{Theorem}
 \label{th1}
Suppose conditions   ${\bf H_1} -{\bf H_2}$ are true.

Then
$$
\zeta_{n}\Longrightarrow \zeta,
$$
where  $\zeta$ is a Poisson point process on $[0,1]\times\mathbb{R}_+$  with   intensity measure  $\mu\times\gamma.$
\end{Theorem}
As has been already noted, this theorem is a generalization of Proposition 3.21  in \cite{r}, which corresponds to the case $c_{ni}=1,\;  \forall\; n,i;\;\;X_{n,i}$ are i.i.d. random elements.

We will prove this theorem in Section \ref{secproof}, and now consider some corollaries.

\section{Argmaxima, maxima and ladder processes}
\renewcommand{\theequation}{\thesection.\arabic{equation}}
\setcounter{equation}{0}

Let $\mathbb{K}$ be a space of locally finite configurations in  $[0,1]\times\mathbb{R}_+.$  In our case, this means that for any    $t>0$,  each configuration has a finite number of points in  $[0,1]\times [t,\infty).$

For each configuration $\varkappa\in \mathbb{K}$ , we define a locally finite measure
$$
\tau := \sum_{x\in \varkappa}\delta_{\{x\}}.
$$
$Vague$-convergence of such measures generates a metric topology in      $\mathbb{K}$ which makes   $\mathbb{K}$   a complete separable space; details may be found, for example, in     \cite{r}.

Consider functionals $A:\mathbb{K}\rightarrow [0,1],\;\;\;M:\mathbb{K}\rightarrow \mathbb{R}_+$ such that
$$
A(\varkappa) = {\mathrm{argmax}}\{x\;| \;(t,x)\in \varkappa \};
$$
$$
M(\varkappa) = \max\{x\;| \;(t,x)\in \varkappa \},
$$
and a map     $L:\mathbb{K}\rightarrow \mathbb{D}$ from  $\mathbb{K}$   to a Skorokhod space  $\mathbb{D}:=\mathbb{D}[0,1]$ for which
 $$
L(\k)(t)=\max\{x\;| \;(s,x)\in \varkappa,\; s\geq t \},\;\;\;t\in [0,1.]
 $$

Since distribution  $\mathcal{P}_\zeta$    of our limiting point process $\zeta$ is defined on  $\mathbb{K},$
and in virtue of the convergence type of configurations in $\mathbb{K},$  (see \cite{r}, Proposition 3.13), functionals  $A,\;M$  and map  $L$ will be almost everywhere (with respect to $\mathcal{P}_\zeta$) continuous. Hence, from Theorem \ref{th1}, we obtain

\begin{Corollary}
	Under the conditions of Theorem \ref{th1},
$$
{\mathrm{argmax}}_{j\leq n}\{X_{n,j}\}  \Longrightarrow  A(\zeta),
$$	
$$
\max_{j\leq n}\{X_{n,j}\}  \Longrightarrow  M(\zeta),
$$	
	$$
	L(\zeta_n)  \Longrightarrow  L(\zeta).
	$$
\end{Corollary}

Next, we find the distribution of variables    $A(\zeta)$   and  $M(\zeta).$

Regarding     $M(\zeta)$,    it is simple: for   $x>0$

$$
\mathbb{P}\{M(\zeta) \leq x\}  = \mathbb{P}\{\zeta([x,\infty]=0)\} = \exp\{-\gamma([x,\infty])\}.
$$
\begin{Proposition}\label{prop1}
	The distribution of  $A(\zeta)$   coincides with $\mu.$
\end{Proposition}

\textbf{Proof.}
	Suppose so far that measure $\gamma$    is finite, set $m=\gamma(\mathbb{R}_+)$, and $\gamma_1 = \frac{1}{m}\gamma.	$

Then, as well known, process  $\zeta$   is equal, in distribution, to a process
$$\pi =: \sum_{j=1}^{\tau} \delta_{\{(Y_j,Z_j)\}},
$$
where r.v.   $\tau$    and sequences   $(Y_j),\;(Z_j)$       are mutually independent, and
\begin{itemize}
  \item the distribution of $\tau $ is Poisson with parameter $m;$
  \item  r.v.'s $Y_j$   are independent, take on values from $[0,1]$ and its distribution equals $\mu;$
  \item r.v.'s   $Z_j$  are independent, non-negative, and its distribution equals  $\gamma_1.$
\end{itemize}

Clearly, the conditional distribution of of   $A(\pi)$  given  $\tau$  and all $Z_j$'s is $\mu$. Then the unconditional distribution of   $A(\pi)$     is   $\mu $ either. Hence, the distribution of  $A(\zeta)$    is also $\mu.$

Let us consider the general case.  Let $\zeta_n,$ be the restriction of  $\zeta$    in  $[0,1]\times [\frac{1}{n},\infty)$.  This is a Poisson point process with intensity measure   $\mu\times\gamma_n,$ where   $\gamma_n$   is the restriction of $\gamma$ on  $[\frac{1}{n},\infty).$

Since $\gamma$  is a Radon measure, for any  $n$, measure $\gamma_n$   is finite.

We know that the distribution of  $A(\zeta_n)$    equals   $\mu$.

 Since with probability one, starting from some $n$, we have  $A(\zeta_n) = A(\zeta),$
 \[
  A(\zeta_n) \Longrightarrow A(\zeta).
  \]
  Therefore the distribution of    $  A(\zeta)$ equals   $\mu.$
  \,\,\,  $\blacksquare$

\vspace{20pt}

\vspace{20pt}

\section{On a connection with paper \cite{d-r}}

\renewcommand{\theequation}{\thesection.\arabic{equation}}
\setcounter{equation}{0}
In this section, we clarify the connection  between the conditions ${\bf H_1}, {\bf H_2}$  of this paper and those of  \cite{d-r}. We will show that under a minor additional condition the integral limit theorem from \cite{d-r} (Theorem 2) may be easily derived from Theorem \ref{th1} of the present paper.

In \cite{d-r},  we considered the following scheme.

First, we defined a sequence  $X_{1},X_{2},...$ of   positive  and  independent r.v.'s,  set $F_{i}(x)=P(X_{i}\leq x)$ and supposed $F(0)=0, F(x)>0$ for all $x>0$.

Next, for  $x>0$, we set
\[
\nu_i(x)=-\ln F_i(x) \,\,\,\,
\]
and $\nu_i(0)=\infty$.

So, for all $i$,
\begin{eqnarray}
	&& F_i(x)=\exp\{-\nu_i(x)\},\label{nu0}   \\
	&& \nu_i(x) \,\,\,\text{ is non-increasing},     \,\,\, \nu_i(0)=\infty, \,\,\, \nu_i(\infty)=0.  \label{nu1}
\end{eqnarray}

The asymptotic behavior of $\nu_i(x) $ as $x\to\infty$ is equivalent to that of $1-F_i(x)$.\vspace{.1in}

Below, we assume all $\nu_i(x)$'s to be strictly decreasing and continuous for $x>0$.

Above this, we impose a condition from \cite{d-r}:

 ${\bf H_3}:$  {\it
  \begin{equation}\label{cond1}
	\nu_i(x)=c_i r(x) (1+\delta_i(x)),
\end{equation}
where $r(x)$ is monotone, all  $\delta_i(x)$ are continuous,  uniformly in $i$
\begin{equation}\label{cond1+}
	\delta_i(x)\to 0 \,\,\, \text{as}\,\,\, x\to\infty,
\end{equation}
and for  positive constants $M<\infty$ and $m<1$, and for all $i$ and $x$,
\begin{equation}\label{cond1++}
	-m\leq \delta_i(x)\leq M.
\end{equation}
}

 It is straightforward to verify that, if $g(x) $ is continuous strictly increasing function, $g(0)=0$ and $g(\infty)=\infty$, then the sequence of r.v.'s    \{ $(\widetilde{X}_j)\}= \{g(X_j)\}$  satisfies condition   ${\bf H_3}$  with corresponding parameters
 \[
\tilde{r}(x)=r(g^{-1}(x)),\;\;\;\tilde{\delta}_j(x)=\delta_j(g^{-1}(x).
\]

On the other hand, clearly,
\[
{\mathrm{argmax}}_{j\leq n}\{X_{j}\} =  {\mathrm{argmax}}_{j\leq n}\{\widetilde{X}_{j}\}.
\]

This remark shows that in the original problem about the distribution of \\
${\mathrm{argmax}}_{j\leq n}\{X_{j}\} $ we can assume that $\tilde{r}$
is a predetermined function. To do this, it suffices to take $g$ so that
$\tilde{r}=r(g^{-1}(x)), $ which is equivalent to the equality\\ $g(x)=\tilde{r}^{-1}(r(x)).$

In what follows we will take  $r(x) = x^{-\alpha }$ for some $\alpha >0,$    and from the sequence  $(X_j)$    let's move on to the triangular array
$$
\{X_{n,j}\},\;\;\; X_{n,j}= d_n^{-\frac{1}{\alpha}}X_j,
$$
where       $d_n= \sum_{j=1}^n c_{j}. $

Let  $F_{n,j}$ and   $\mathcal{P}_{n,j}$  be respectivly the distribution function and distribution of    $X_{n,j}.$  Then, due to
${\bf H_3},$
 \begin{equation}\label{h3.1}
	\frac{d_n}{c_j}[1-F_{n,j}(x)]= x^{-\alpha}(1+\theta_{n,j}(x)),	
\end{equation}
where for sufficiently large $n,$ for $j\leq n$ and for $x\geq a$
 \begin{equation}\label{h3.2}
	|\theta_{n,j}(x)| \leq  C_a\frac{M_n}{d_n},
\end{equation}
and

$ M_n = \max_{j\leq n}\{c_j\}.$

      Set of conditions ${\bf H1} $   now looks like:

$ {\bf H1}$:  {\it As $n\rightarrow \infty$
\begin{equation}\label{c1}
	d_n:= \sum_{j=1}^n c_{j}  \rightarrow \infty;\;\;\;\;\;\;
	\frac{\max_{j\leq n}\{c_{j}\}}{d_n} \rightarrow 0,
\end{equation}
and mesures on $[0,1],$
$$
\mu_n:= \sum_{j=1}^n 	\frac{c_{j}}{d_n}\delta_{\{\frac{j}{n}\}},
$$
weakly converge to some measure $\mu$:
\begin{equation}\label{h1}
	\mu_n{\Longrightarrow} \mu.
\end{equation}
}
\vspace{20pt}

\begin{Theorem}
 \label{th3}
Let us assume that the conditions $ {\bf H1,\;\bf H3}$ are fulfild.
Let us also assume that the distribution functions $F_j$ are strictly monotonic. 
Then 
\begin{equation}\label{t2}
{\mathrm{argmax}}_{j\leq n}\{X_{n,j}\} =
{\mathrm{argmax}}_{j\leq n}\{d_n^{-\frac{1}{\alpha}}X_j\}
\Longrightarrow  \mu.
\end{equation}	
\end{Theorem}

\textbf{Proof. }
To prove it, it is enough to note that from the relations ({\ref{h3.1}}),
({\ref{h3.2}}) follows (\ref{h2.2}). Indeed, ({\ref{h3.1}})
means that

$$
\frac{d_n}{c_j}\mathcal{P}_{n,j} \stackrel{vague}{\longrightarrow} \gamma,
$$
where         $\gamma([t,\infty)) = t^{\alpha},\;\; t>0.$

By virtue of ({\ref{h3.2}}), this convergence is uniform in the sense of ({\ref{h2.2}}).

Thus, one can apply Th.{\ref{th1}} and, as a consequence, obtain
({\ref{t2}}) . $\blacksquare$

\section{Proof of Theorem \ref{th1}}\label{secproof}
As is known (see, for example,  \cite{r}), to prove a convergence of point processes, it suffices to establish the convergence of the corresponding Laplace functionals (L.f.).

Let  $f\in \mathbb{C}_{K}^+([0,1]\times \mathbb{R}_+). $ For  $\zeta_n$  the   L.f.\
\begin{eqnarray}
\Psi_n(f)&&= \mathbf{E}\exp\{-\zeta_{n}(f)\}= \mathbf{E}
\exp\left\{-\sum_{j\leq n}f\left(\frac{j}{n},\;X_{n,j}\right)\right\}\\
&&=\prod_{j\leq n}\left\lbrace 1-\int_{\mathbb{R}_+} \left(1-e^{-f(\frac{j}{n},x)} \right)\mathcal{P}_{n,j}(dx)\right\rbrace .
\end{eqnarray}
On the other hand, for $\zeta$  , as is known
\begin{equation}
	\Psi(f)= \exp {\left\lbrace  - \int_{0}^{1}\int_{\mathbb{R}_+}\left(1-e^{-f(t,x)} \right)
	\tau(dt,dx)\right\rbrace }.
\end{equation}
We have
$$
-\ln \Psi_n(f) = - \sum_{j\leq n}\ln\left( 1-\int_{\mathbb{R}_+}\left(1-e^{-f(\frac{j}{n},x)} )\right)\mathcal{P}_{n,j}(dx)\right).
$$

Setting $Q_{n,j}= \frac{d_n}{c_{n,j}}\mathcal{P}_{n,j},$       we get
\begin{eqnarray*}
\Sigma_n:&&= \sum_{j\leq n} \int_{\mathbb{R}_+}\left(1-e^{-f(\frac{j}{n},x)} \right)\mathcal{P}_{n,j}(dx)\\
&&
=\sum_{j\leq n}\int_{\mathbb{R}_+}\left(1-e^{-f(\frac{j}{n},x)} \right)\frac{c_{n,j}}{d_n}\cdot\frac{d_n}{c_{n,j}}\mathcal{P}_{n,j}(dx)\\
&&=\int_{0}^{1}\int_{\mathbb{R}_+}\left(1-e^{-f(t,x)} \right)\tau_n(dt,dx),
\end{eqnarray*}
where measure   $\tau_n$ is defined on  $[0,1]\times \mathbb{R}_+,$ and
$$
\tau_n(A) = \sum_{j\leq n}\mu_n(\frac{j}{n})Q_n\{A\cap({\frac{j}{n}\times}\mathbb{R}_+)\},\;\;\;
A\subset [0,1]\times \mathbb{R}_+.
$$
\begin{Lemma}\label{L}
	As $n\rightarrow \infty$
	\begin{equation}
\tau_n	\stackrel{vague}{\longrightarrow} \mu\times\gamma.	
	\end{equation}
\end{Lemma}
{\bf Proof of Lemma.}
Let $h \in\mathbb{C}_{K}^+([0,1]\times \mathbb{R}_+).$ Then    $\mathrm {supp}\{h\}$  belongs to  $[0,1]\times[a,\infty)$ for some  $a>0$.

In virtue of (\ref{h2.2}), as  $n \rightarrow \infty$,
$$
 \left|  \int h d\tau - \int h d\mu_nd\gamma\right|
 $$
 \begin{equation}\label{l1}
	= \left| \sum_{j\leq n}\left[ \int_{\mathbb{R}_+}h(\frac{j}{n},x)Q_{n,j}(dx) -
\int_{\mathbb{R}_+}h(\frac{j}{n},x)\gamma(dx) \right]\frac{c_{n,j}}{d_n}\right|  
\leq \Delta_n(a) \rightarrow 0.
\end{equation}

We have also
\begin{equation}\label{l2}
	\left|  \int h d\mu_nd\gamma - \int h d\mu d\gamma \right| \leq
\int_{\mathbb{R}_+}H_nd\gamma,
	\end{equation}
where
$$
H_n=		\left|  \int_0^1 h d\mu_n - \int_0^1 h d\mu \right|.
	$$
As  $\mu_n \Rightarrow \mu,$    then  $\forall x>0\;\;\;H_n(x )\rightarrow 0, \;\;n \rightarrow \infty.$
Moreover, there is $C>0$ such that
  $H_n(x )\leq C\mathbf{1}_{[a,\infty)}(x).$
  That's why
  from (\ref{l1}) and (\ref{l2} ) the proof of the lemma follows.  $\blacksquare$

By Lemma \ref{L}
\begin{equation}\label{pr0}
\Sigma_n \longrightarrow  \int_{0}^{1}\int_{\mathbb{R}_+}\left(1-e^{-f(t,x)} )\right)
\tau(dt,dx),
\end{equation}
where    $\tau= \mu\times\gamma.$

Let us show that $\Sigma_n$ approaches   $\;\;-\ln \Psi_n(f).$

Note that from (\ref{h2.2}) it follows that for all $f$  from $ \mathbb{C}_{K}^+( \mathbb{R}_+)$ having support lying in $[ a,\infty),$ and for all sufficiently large $n$

\begin{equation} \label{pr1}
\max_{j\leq n}\left\lbrace  \frac{d_n}{c_{n,j}}\int_{\mathbb{R}_+} f d\mathcal{P}_{n,j}\right\rbrace \;\;\;
\leq  \;\;\; 2\int_{\mathbb{R}_+} f d\gamma  \;\;\;\leq  \;\;\;\|f\|_{\infty} \gamma([a,\infty)).
\end{equation}
As for   $|t|\leq 1/2$
$$
\ln(1+t)=t(1+\varepsilon(t)),\;\;\;|\varepsilon(t)| \leq |t|,
$$
then for all sufficiently large  $n$
\begin{multline}\label{pr2}
	\left| -\ln \Psi_n(f) - \Sigma_n\right|\\
	 =  \left| - \sum_{j\leq n}\ln\left( 1-\int_{\mathbb{R}_+}\left(1-e^{-f(\frac{j}{n},x)} )\right)\mathcal{P}_{n,j}(dx)\right) -
\sum_{j\leq n} \int_{\mathbb{R}_+}\left(1-e^{-f(\frac{j}{n},x)} )\right)\mathcal{P}_{n,j}(dx) \right|\\
\leq  \;\;\; \sum_{j\leq n} \left( \int_{\mathbb{R}_+}\left(1-e^{-f(\frac{j}{n},x)} )\right)\mathcal{P}_{n,j}(dx)\right)^2 \\
\leq  \;\;\;\max_{j\leq n}\left\lbrace \int_{\mathbb{R}_+}\left(1-e^{-f(\frac{j}{n},x)} )\right)\mathcal{P}_{n,j}(dx)   \right\rbrace \cdot
\sum_{j\leq n} \int_{\mathbb{R}_+}\left(1-e^{-f(\frac{j}{n},x)} )\right)\mathcal{P}_{n,j}(dx).
\end{multline}
Due to (\ref{h2.1})  for all sufficiently large  $n$
$$
\int_{\mathbb{R}_+}\left(1-e^{-f(\frac{j}{n},x)} \right)\mathcal{P}_{n,j}(dx)   \;\;\;  \leq \;\;\; \mathcal{P}_{n,j}([a,\infty)) \;\;\;\leq \;\;\; 2a^{-\alpha }\frac{M_n}{d_n} \rightarrow 0,\;\;\;n
\rightarrow \infty;
$$
here $ M_n = \max_{j\leq n}\{c_{n,j}\}.$

Therefore,  taking (\ref{pr1}) into account, we obtain from  (\ref{pr2}) that there is a constant $C,$ depending only on $f,$ such that

$$
\left| -\ln \Psi_n(f) - \Sigma_n\right|\leq C\frac{M_n}{d_n}\Sigma_n \rightarrow 0,\;\;\;n
\rightarrow \infty.
$$
By virtue of (\ref{pr0}) we finally get
$$
\Psi_n(f)  \rightarrow  \Psi(f),
$$
which proves the theorem.  $\blacksquare$

\end{document}